\documentclass[12pt]{article}

\usepackage{brunnian}

\usepackage[english,francais,catalan]{babel}
\usepackage[T1]{fontenc}
\usepackage[utf8]{inputenc}
\usepackage{geometry}
\usepackage{amsmath}
\usepackage{amsthm}
\usepackage{thmtools}
\usepackage{amscd}
\usepackage{amsfonts}
\usepackage{amssymb}
\usepackage[all,cmtip]{xy}

\usepackage{hyperref}
  \hypersetup{
      colorlinks,
      citecolor=black,
      filecolor=black,
      linkcolor=black,
      urlcolor=black
  }
\usepackage{nameref,cleveref}
\usepackage{float}
\usepackage{epsfig}
\usepackage{latexsym}
\usepackage{verbatim}
\usepackage{indentfirst}
\usepackage{graphicx}
\usepackage{caption}
\usepackage{subcaption}
\usepackage{cancel}
\usepackage{xspace}
\usepackage{times}
\usepackage{upgreek}
\usepackage{catchfilebetweentags}
\usepackage{wasysym}
\usepackage[section]{placeins}
\usepackage{array}
\usepackage{lipsum}
\usepackage{alltt}

  \def\C{\mathbb{C}}
  \def\N{\mathbb{N}}

  \def\S{\mathbb{S}}
  \def\Z{\mathbb{Z}}

  \def\Aa{\mathcal{A}}

  \def\Ll{\mathcal{L}}
  
\def\ol{\overline}
\def\ul{\underline}

\def\>{\rangle}
\def\<{\langle}
\def\={\neq}

\def\l{\ell}

\def\dim{\mathrm{dim}}

\def\ker{\mathrm{Ker}}
\def\tr{\mathrm{tr}}
\def\hom{\mathrm{Hom}}

\def\ab{\mathrm{ab}}
\def\red{\mathrm{red}}

\def\mf{\mathfrak{m}}
\def\lf{\mathfrak{l}}

\def\-{\!\!-\!\!}

\def\edge#1#2#3{\xymatrix@C=8pt{#1\ar@{-}[r]^-{#2} & #3}}

\def\Tarrow#1#2{\xymatrix@C=14pt{#1\ar[r]^-{#2} &}}

\def\contextEdge#1#2#3#4#5{\xymatrix@C=8pt{#1\supset #2\ar@{-}[r]^-{#3} & #4\subset #5}}

\DeclareMathOperator{\lk}{\mathrm{lk}}

\renewcommand{\geq}{\geqslant}
\renewcommand{\leq}{\leqslant}

\declaretheorem{theorem}
\declaretheorem[sibling=theorem]{lemma}
\declaretheorem[sibling=theorem]{proposition}

\declaretheorem[sibling=theorem,style=definition]{definition}

\declaretheorem[style=remark,numberwithin=section]{remark}

\title{Eigenvalue varieties of Brunnian links}
\author{François Malabre}

\begin{document}
\selectlanguage{english}

\maketitle


\begin{abstract}

  In this article, it is proved that the eigenvalue variety of the exterior of a nontrivial, non-Hopf, Brunnian link in $\S^{3}$ contains a nontrivial component of maximal dimension.
  Eigenvalue varieties were first introduced to generalize the $A$-polynomial of knots in $\S^{3}$ to manifolds with nonconnected toric boundary.
  The result presented here generalizes, for Brunnian links, the nontriviality of the $A$-polynomial of nontrivial knots in $\S^{3}$.

\end{abstract}



The $A$-polynomial of a knot in $\S^{3}$ is a two-variable polynomial constructed from the $\sl$-character variety of the knot exterior.
Let $K$ be a knot in $\S^{3}$ and let $\pi_{1}K$ denote the fundamental group of the exterior of $K$;
the peripheral subgroup $\Z^{2}$ is generated by a meridian $\mu$ and a longitude $\lambda$, and the zero-set of the $A$-polynomial $A_{K}$ is the locus of eigenvalues for a common eigenvector of $\rho(\mu)$ and $\rho(\lambda)$ of representations $\rho$ from $\pi_{1}K$ to $\sl$.
It was first introduced by Cooper, Culler, Gillet, Long and Shalen in \cite{Cooper:1994gp}, where it is also proved that the $A$-polynomial of any knot contains the $A$-polynomial of the unknot as a factor.
The $A$-polynomial of a knot is said to be \emph{nontrivial} if it contains other factors and it was also proved, in the same article \cite{Cooper:1994gp}, that hyperbolic knots and nontrivial torus knots always have a nontrivial $A$-polynomial.
This was later established in full generality for all nontrivial knots by Dunfield and Garoufalidis in \cite{Dunfield:2004ei}, and independently by Boyer and Zhang in \cite{Boyer:2005hu};
both proofs use a theorem by Kronheimer and Mrowka in \cite{Kronheimer:2004jr} on Dehn fillings on knots and representations in $\su$.

The notion of $A$-polynomial can be generalized to any $3$-manifold $M$ with connected toric boundary by specifying a \emph{peripheral system} (generators of $\pi_{1}\partial M \hookrightarrow \pi_{1}M$).
Stimulated by the work of Lash in \cite{lash1993boundary}, it was then extended to manifolds with non-connected boundary by Tillmann.
In his PhD thesis \cite{Tillmann:2002tz} and the subsequent article \cite{Tillmann:2005ck}, Tillmann presented the \emph{eigenvalue variety} $\evar{M}$ associated to a $3$-manifold $M$ with toric boundary.
If the boundary of $M$ consists of $n$ tori, the associated eigenvalue variety $\evar{M}$ is an algebraic subspace of $\C^{2n}$ corresponding to the closure of peripheral eigenvalues taken by representations (or equivalently, characters) of $\pi_{1}M$ in $\sl$.
Under these assumptions, Tillmann proved in \cite{Tillmann:2005ck} that the dimension of any component of $\evar{M}$ is at most $n$.

In the same way as any $A$-polynomial is divisible by the $A$-polynomial of the unknot, any eigenvalue variety $\evar{M}$ contains components $\evarGen{\red}{}{M}$ corresponding to reducible characters.
Components of $\evarGen{\red}{}{M}$ have maximal dimension and any other component of $\evar{M}$ with maximal dimension is called a \emph{nontrivially maximal} component of $\evar{M}$.

If $M$ is hyperbolic, its character variety contains a family of distinguished components $Y_{1},\dots, Y_{k}$ called the \emph{geometric components}, each one containing the character of a discrete faithful representation.
Using Thurston's results of \cite{Thurston:tj}, Tillmann proved that each geometric component produces a nontrivially maximal component in $\evar{M}$, generalizing the result of \cite{Cooper:1994gp} on hyperbolic knots.
However, for which $3$-manifolds $M$ does $\evar{M}$ contain a nontrivially maximal component, or merely whether this is true or not for nontrivial exteriors of links in $\S^{3}$, remain open questions.

In this article, we answer this matter for a family of links in $\S^{3}$, the \emph{Brunnian links}.
A link in $\S^{3}$ is called \emph{Brunnian} if any of its proper sublinks is trivial and we prove the following:
\begin{restatable}{theorem}{brunnMax}
\label{th:brunnian_links_maximal_evar}
  The eigenvalue variety of any nontrivial non-Hopf Brunnian link contains a nontrivially maximal component.
\end{restatable}

The defining property of Brunnian links makes them stable under $1/q$-Dehn fillings, which permits to apply Kronheimer-Mrowka's Theorem to produce irreducible characters in a similar fashion as in \cite{Dunfield:2004ei} and \cite{Boyer:2005hu}.
Then, an induction on the number of components of the links produces nontrivially maximal components of their eigenvalue varieties.

This article is divided into two sections;
first we recall the construction of the eigenvalue variety $\evar{L}$ for a link $L$ in $\S^{3}$, its defining ideal $\aid{L}$ and some  of its properties, as presented in \cite{Tillmann:2005ck}, to introduce notation for the following section.
Then, we study the family of Brunnian links in $\S^{3}$ and prove the main result of this article.



\section{Eigenvalue varieties of links in $\S^3$} 
\label{sec:eigenvalue_varieties_of_links_in_s_3}

  First, we briefly review the notion of \emph{eigenvalue variety} associated to a link in $\S^{3}$;
  this was first introduced by Tillmann in \cite{Tillmann:2002tz} and we reproduce the construction here (with a slightly different vocabulary) in order to set the notation for the next section.

  \subsection{Character varieties} 
  \label{sub:character_varieties}

Let $\pi$ be a finitely generated group;
the \emph{$\sl$-representation variety} of $\pi$ is the algebraic affine set $\hom(\pi,\sl)$ and is denoted by $R(\pi)$.
The algebraic Lie group $\sl$ acts on $R(\pi)$ by conjugation and the algebraic quotient under this action is the \emph{$\sl$-character variety} of $\pi$, denoted by $X(\pi)$.
The ring of regular functions on the character variety, $\C[X(\pi)]$, is equal to the subring $\C[R(\pi)]^{\sl}$ of invariant functions.
Dually, the inclusion $\C[X(\pi)]\hookrightarrow \C[R(\pi)]$ induces a natural algebraic epimorphism $t\colon R(\pi)\to X(\pi)$ and any regular function on $R(\pi)$ factors through $t$ if and only if it is invariant under the conjugation action of $\sl$.
In particular for any $\gamma$ in $\pi$, the function $\tau_{\gamma}\colon R(\pi)\to \C$ mapping $\rho \mapsto \tr~\rho(\gamma)$ defines a regular function $I_{\gamma}$ on $X(\pi)$ called the \emph{trace function at $\gamma$};
the trace functions finitely generate the ring $\C[X(\pi)]$ (see \cite{Culler:1983gka} for example).
Representation and character varieties are contravariant functors:
any group morphism $\pi\to\pi'$ induces regular maps according to the following commutative diagram:
\[
  \xymatrix{
    R(\pi') \ar[r]\ar[d]_-{t} & R(\pi)\ar[d]^-{t} \\
    X(\pi') \ar[r] & X(\pi).
  }
\]

In case the group $\pi$ is the fundamental group of a manifold $M$ (resp. the exterior of a link $L$ in $\S^{3}$), the representation and character varieties will be denoted by $R(M)$ and $X(M)$ (resp. $R(L)$ and $X(L)$).



  \subsection{Abelian characters} 
  \label{sub:abelian_characters}

Any group $\pi$ has an abelianization $\pi^{\ab}$ and a canonical projection $\pi\to\pi^{\ab}$ which induces regular maps
\[
  \xymatrix{
    R(\pi^{\ab}) \ar[r]\ar[d]_-{t} & R(\pi)\ar[d]^-{t} \\
    X(\pi^{\ab}) \ar[r] & X(\pi).
  }
\]
The image of $R(\pi^{\ab})$ in $R(\pi)$ is precisely the closed set $R^{\ab}(\pi)$ of abelian representations of $\pi$ and the image of $X(\pi^{\ab})$ is a closed subset of $X(\pi)$ called the set of \emph{abelian characters} of $\pi$ and denoted by $X^{\ab}(\pi)$.

\begin{remark}
\label{rem:red_chars_are_ab_chars}
  In $\sl$, characters of reducible representations are characters of abelian representations.
  If $R^{\red}(\pi)$ is the closed set of reducible representations and $X^{\red}(\pi)$ is its image in $X(\pi)$, then $X^{\red}(\pi) = X^{\ab}(\pi)$.
\end{remark}

Let $\Delta$ denote the map from $\C^{\ast}$ to $\sl$ mapping $z\mapsto \DeltaMat{z}$;
by composition, $\Delta$ defines maps
\[
  \xymatrix{
  \hom(\pi,\C^{\ast}) \ar[r]^-{\Delta_{\ast}}\ar[dr]_-{d} & R^{\ab}(\pi)\ar[d]^-{t}\\
   & X^{\ab}(\pi).
  }
\]
The map $d$ is $2:1$ onto $X^{\ab}(\pi)$, invariant under inversion in $\hom(\pi,\C^{\ast})$;
for any $\varphi$ in $\hom(\pi,\C^{\ast})$ and $\gamma$ in $\pi$,
\[
  I_{\gamma}\circ d(\varphi) = \traceDiag{\varphi(\gamma)}.
\]



  \subsection{Eigenvalue varieties} 
  \label{sub:eigenvalue_varieties}

Let $L$ be a link in $\S^{3}$, let $|L|$ denote its number of components and let $\pi_{1}L$ be the fundamental group of its exterior;
the boundary of the exterior of $L$ is a disjoint union of $|L|$ tori $T_{K}$, one for each component $K$ of the link $L$.
Each inclusion $\pi_{1}T_{K}\hookrightarrow \pi_{1}L$ induces a regular map $r_{K}:X(L)\to X(T_{K})$.
Since $\pi_{1}T_{K}$ is abelian, $X(T_{K})=X^{\ab}(T_{K})$ and denoting $\hom(\pi_{1}T_{K},\C^{\ast})$ by $E(T_{K})$ we obtain the following diagram:
\[
  \xymatrix{
    & \prod_{K\subset L}E(T_{K})\ar[d]^-{d}\\
    X(L) \ar[r]_-{r} & \prod_{K\subset L}X(T_{K}).
  }
\]
Following Tillmann \cite{Tillmann:2002tz,Tillmann:2005ck}, the \emph{eigenvalue variety} of $L$ is defined as the Zariski closure of the preimage by $d$ of the image of $r$:
\[
  \evar{L}=\ol{d^{-1}(r(X(L)))}.
\]

Dually, there are ring-maps
\[
  \xymatrix{
    & \otimes_{K\subset L}\C[E(T_{K})]\\
    \C[X(L)] & \otimes_{K\subset L}\C[X(T_{K})] \ar[l]^-{r^{\ast}}\ar[u]_-{d^{\ast}}
  }
\]
and the defining ideal $\aid{L}$ of $\evar{L}$ is called the \emph{$\Aa$-ideal} of $L$ and is the radical of the image by $d^{\ast}$ of the kernel of $r^{\ast}$:
\[
  \aid{L}=\sqrt{d^{\ast}(\ker~r^{\ast})}.
\]

Each torus $T_{K}$ is equipped with a \emph{standard peripheral system} $(\mu_{K},\lambda_{K})$ of meridian and longitude of each component.
This produces canonical coordinates $(m_{K},\l_{K})$ in $\C^{\ast}\times \C^{\ast}$ for $E(T_{K})$, and $\evar{L}$ is naturally a subset of $(\C^{\ast})^{2|L|}$;
dually, $\C[E(T_{K})]$ is isomorphic to $\C[\mf_{K}^{\pm 1},\lf_{K}^{\pm 1}]$ and $\aid{L}$ is an ideal of $\C[\mf^{\pm},\lf^{\pm}] = \otimes_{K\subset L}\C[\mf_{K}^{\pm 1},\lf_{K}^{\pm 1}]$.

\begin{proposition}
\label{prop:A^red}
  Let $\evarGen{\ab}{}{L}$ denote the part of $\evar{L}$ corresponding to abelian characters and $\aidGen{\ab}{}{L}$ the corresponding defining ideal;
  $\evarGen{\ab}{}{L}$ is a union of copies of $(\C^{\ast})^{|L|}$ and $\aidGen{\ab}{}{L}$ is given in $\C[\mf^{\pm},\lf^{\pm}]$ by
  \[
    \aidGen{\ab}{}{L} = \left\langle \lf_{K} - \prod_{K'\neq K}\mf_{K'}^{\pm\lk(K,K')} \right\rangle
  \]
  where $\lk(K,K')$ denotes the linking number of the components $K$ and $K'$.
\end{proposition}

\begin{proof}
  The meridians form a basis of the homology group of the link exterior and each longitude is given by the linking numbers
  \[
    \lambda_{K}=\sum_{K'\neq K}\lk(K,K')\mu_{K'}.
  \]
  Therefore, any morphism from $\pi_{1}L$ to $\C^{\ast}$ is determined by the images of the meridians and for any $\varphi$ in $\hom(\pi_{1}L,\C^{\ast})$ and each longitude $\lambda_{K}$,
  \[
    \varphi(\lambda_{K}) = \prod_{K\neq K'}\varphi(\mu_{K'})^{\lk(K,K')}.
  \]

  By the invariance under inversion, any point $(m_{K},\l_{K})_{K\subset L}$ of $\evarGen{\ab}{}{L}$ satisfies then
  \[
    \l_{K} = \prod_{K\neq K'}{m_{K}}^{\pm\lk(K,K')}.
  \]

  Conversely, for any $\xi = (m_{K},\l_{K})_{K\subset L}$ satisfying these equations, there exists $\varphi$ in $\hom(\pi_{1}L,\C^{\ast})$ such that $d(\xi) = r(\Delta_{\ast}\varphi)$ so $\aidGen{\ab}{}{L}$ is given by
  \[
    \aidGen{\ab}{}{L} = \left\langle \lf_{K} - \prod_{K'\neq K}\mf_{K'}^{\pm\lk(K,K')} \right\rangle.
  \]
\end{proof}

\begin{remark}
\label{rem:A^red_generalizes_l-1}
  For links with one component (knots), the $\aidC$-ideal is generated by the $A$-polynomial of the knot and $A^{\ab}$ is the $\lf-1$ factor corresponding to abelian characters.
\end{remark}

By the defining equations of $\aidGen{\ab}{}{L}$, $\evarGen{\ab}{}{L}$ always has dimension $|L|$.
As a matter of fact, by Tillmann \cite{Tillmann:2002tz,Tillmann:2005ck}, any component of $\evar{L}$ has dimension at most $|L|$, which leads to the following definition:

\begin{definition}
\label{def:nontrivially_maximal}
  A component of $\evar{L}$ is called \emph{nontrivially maximal} if it has dimension $|L|$ and is not contained in $\evarGen{\ab}{}{L}$.
\end{definition}

Using Thurston's results on hyperbolic manifolds, Tillmann showed the following:
\begin{theorem}[Tillmann \cite{Tillmann:2005ck}]
\label{th:hyp_link_max_evar}
  If $L$ is a hyperbolic link in $\S^{3}$ then any geometric component of the character variety produces a nontrivially maximal component in the eigenvalue variety of $L$.
\end{theorem}

Besides these cases, it is not known whether the eigenvalue variety of all (nontrivial) links admits a maximal nontrivial component.
For knots, this is equivalent to the nontriviality of the $A$-polynomial (besides the $\lf-1$ factor) and was proven idependently by Dunfield-Garoufalidis in \cite{Dunfield:2004ei} and Boyer-Zhang in \cite{Boyer:2005hu}.
In the next section, we answer this matter for Brunnian links in $\S^{3}$.






\section{Characters of Brunnian links} 
\label{sec:characters_of_brunnian_links}

  In this section we prove \Cref{th:brunnian_links_maximal_evar}.
  First, we recall some basic facts on $1/q$-Dehn fillings on links in $\S^{3}$;
  then, we present Brunnian links and, after having studied their stability under these Dehn fillings, we use Kronheimer-Mrowka's Theorem to create families of characters of Brunnian links exteriors.
  Finally, we prove that these characters span a nontrivially maximal component in the eigenvalue varieties of nontrivial, non-Hopf, Brunnian links.

  \subsection{Dehn fillings} 
  \label{sub:dehn_fillings}

Any $1/q$-Dehn filling on the unknot in $\S^{3}$ produces $\S^{3}$ again;
therefore, the $1/q$-Dehn filling over an unknotted component of a link in $\S^{3}$ produces the exterior of another link in $\S^{3}$.

Let $L=K\sqcup L'$ be a link with $K$ an unknotted component of $L$, and let $L_{q}$ denote the link obtained by $1/q$-surgery on $K$ (so, in particular, $L'=L_{0}$).
Any sublink $L''$ of $L_{q}$ is obtained by $1/0$-Dehn filling along the other components.
Because the meridians are unchanged by $1/q$-Dehn fillings, any proper sublink $L''$ of $L_{q}$ is obtained by $1/q$-Dehn filling along $K$ in the sublink $L''\sqcup K$ of $L$.
\begin{remark}
\label{rem:sublink_trivial_dehn}
  With this notation, if $L''\sqcup K$ is trivial in $\S^{3}$, then so is $L''$.
\end{remark}

The meridians are unchanged by $1/q$-Dehn fillings but the longitudes are changed according to the linking numbers.
With the same notation as above, if $(\mu,\lambda)$ is a standard peripheral system for a component $J$ of $L$, then the new longitude $\lambda_{q}$ of $J$ in $L_{q}$ is
\[
  \lambda_{q} = \lambda + q~\lk(K,J)^{2}\mu
\]
and the linking number $\lk_{q}(J,J')$ of any two components $J$ and $J'$ of $L_{q}$ is given by
\[
  \lk_{q}(J,J')=\lk(J,J')-q~\lk(K,J)~\lk(K,J').
\]

A link is called \emph{homologically trivial} if all the linking numbers between components vanish.
By the previous discussion, the link obtained by $1/q$-Dehn fillings on an unknotted component of a homologically trivial link is still homologically trivial and has the same longitudes.

The proof of \Cref{th:brunnian_links_maximal_evar} uses Dehn fillings to produce closed $3$-manifolds which admit irreducible representations;
this will be done by iterating $1/q$-Dehn fillings along the components of the link.
However, even if all the components of a link $L$ in $\S^{3}$ are unknotted, a $1/q$-Dehn filling along a component generally knots the other components, thus making impossible to continue the process while remaining in $\S^{3}$.
In other words, to achieve this goal, we need a family of links $\Ll$ satisfying:
\begin{itemize}
  \item if $L\in\Ll$ has two or more components, each is individually unknotted;
  \item for any $K\sqcup L_{0}$ in $\Ll$, $L_{q}$ is also in $\Ll$.
\end{itemize}
In the next section, we show that the family of \emph{Brunnian links} in $\S^{3}$ satisfies these conditions.
Moreover, nontriviality can be preserved in the process, making it possible to reason by induction on the number of components of the link.



  \subsection{Brunnian links} 
  \label{sub:brunnian_links}

\begin{definition}
\label{def:brunnian_link}
  A link is called \emph{Brunnian} if any of its proper sublinks is trivial.
\end{definition}

\begin{remark}
\label{rem:brunnian_link}
  Any knot is considered Brunnian;
  for links with more components we have:
  \begin{itemize}
    \item The components of a Brunnian link with $2$ components or more are individually unknotted.
    \item Any Brunnian link with $3$ or more components is homologically trivial.
    \item By \Cref{rem:sublink_trivial_dehn}, if $L=K\sqcup L_{0}$ is Brunnian, $L_{q}$ is also Brunnian for any integer $q$.
  \end{itemize}
\end{remark}

Given $L = K\sqcup L_{0}$ Brunnian, we can perform a $1/p$-surgery on a component of $L_{q}$ to obtain another Brunnian link, and so on, until obtaining a knot in $\S^{3}$.
However, any $1/q$-Dehn filling on a component of the Hopf link or the unlink produces the unlink.
Therefore, given a Brunnian link $L=K\sqcup L_{0}$, we need to prevent $L_{q}$ from being the Hopf link or the unlink in order to obtain, \emph{in fine}, a nontrivial knot in $\S^{3}$.

If $L=K\sqcup K'$ is a Brunnian link with two components, this is a special case of Mathieu's Theorem from \cite{Mathieu:1992uv}.
This more general result on knots in a solid torus (links with one unknotted component) asserts that, besides the Hopf link, for any $|q|\geq 2$, any $1/q$-Dehn filling on an unknotted component of a $2$-component link in $\S^{3}$ produces a nontrivial knot.
For our concern, this implies that, for any $|q|\geq 2$, the $1/q$-Dehn filling on any component of a Brunnian, non-Hopf, nontrivial $2$-link may never produce the trivial knot.

On the other hand, if $L$ has three components or more, it is homologically trivial and the work of Mangum-Stanford in \cite{Mangum:2001cc} (Theorem 2 and its proof) ensures that, for any integer $q$ and any homologically trivial Brunnian link $L=K\sqcup L_{0}$, if $L$ is nontrivial, then $L_{q}$ is trivial if and only if $q=0$.
Otherwise, it is a nontrivial, homologically trivial Brunnian link (in particular, it is never the Hopf link).

Therefore, we obtain the following result for the stability of nontrivial non-Hopf Brunnian links under $1/q$-Dehn fillings:
\begin{proposition}
\label{prop:non_trivial_brunnian_links_fillings}
  Let $L=K\sqcup L_{0}$ be a nontrivial, non-Hopf, Brunnian link in $\S^{3}$.
  Then, for any $|q|\geq 2$ the link $L_{q}$ is a Brunnian link in $\S^{3}$, nontrivial and non-Hopf.
\end{proposition}

We will use the stability of nontrivial non-Hopf Brunnian links to apply Kronheimer-Mrowka's theorem on some Dehn fillings of the link exteriors to produce nontrivially maximal components in the eigenvalue varieties.
On the other hand, for the Hopf link and the trivial link, no such component exists:

\begin{proposition}
\label{prop:Hopf_unlink_trivial_evar}
  The eigenvalue varieties of the Hopf link and the trivial link do not admit any nontrivially maximal component.
\end{proposition}

\begin{proof}
  The fundamental group of the exterior of the Hopf link is abelian so all the characters are abelian and $\evarC=\evarC^{\ab}$.

  On the other hand, for the trivial link, all the longitudes are nullhomotopic and are therefore trivialized by any representation so the $\aidC$-ideal is $\<\lf_{K}-1,~K\subset L\>=\aidC^{\ab}$.
\end{proof}



  \subsection{Kronheimer-Mrowka characters} 
  \label{sub:kronheimer_mrowka_characters}

By Kronheimer-Mrowka's Theorem from \cite{Kronheimer:2004jr}, any nontrivial $1/q$-Dehn filling along a nontrivial knot in $\S^{3}$ produces a closed $3$-manifold wich admits an irreducible representation in $\su$.
By \Cref{prop:non_trivial_brunnian_links_fillings}, if $L = K\sqcup L_{0}$ is a nontrivial Brunnian link in $\S^{3}$, $L_{q}$ is nontrivial for any $|q|\geq 2$.
Performing another $1/p$-Dehn filling on a component of $L_{q}$ (in the new standard peripheral system if the link is not homologically trivial) will produce again a nontrivial Brunnian link;
this process may be continued until a nontrivial knot is produced, on which a final $1/k$-Dehn filling may be performed to obtain a closed $3$-manifold wich admits an irreducible representation in $\su$.

For any Brunnian link $L=K_{1}\sqcup\dots\sqcup K_{n}$ in $\S^{3}$, and any $\ul{q}=(q_{1},\dots,q_{k})$ in $\Z^{k}$ for $k\leq n$, we denote by $L(\ul{q})$ the $3$-manifold obtained by performing $1/q_{i}$-Dehn fillings on the components of $L$, where each $1/q_{i}$-Dehn filling is performed in the standard peripheral system given after the Dehn fillings $1/q_{j}$ for $j<i$.
\begin{remark}
\label{rem:std_long_changes}
  As already pointed out, the meridians never change and, since $L$ is assumed Brunnian, longitudes change only if $L$ is a Brunnian link with two components $L=K_{1}\sqcup K_{2}$ with nonzero linking number $\alpha$;
  in that case, denoting by $(\mu_{i},\lambda_{i})_{i=1,2}$ the respective standard peripheral systems, any $1/q_{1}$-Dehn filling on $K_{1}$ changes the longitude $\lambda_{2}$ into $\lambda_{2}+q_{1}\alpha^{2}\mu_{2}$.
  Therefore, a $1/q_{2}$-Dehn filling on $K_{2}$ is performed along the slope
  \[
    (1+q_{1}q_{2}\alpha^{2})\mu_{2} + q_{2}\lambda_{2} \in H_{1}(T_{K_{2}}).
  \]
\end{remark}

\begin{proposition}
\label{prop:kronheimer_mrowka_brunn}
  Let $L=K_{1}\sqcup\dots\sqcup K_{n}$ be a nontrivial Brunnian link in $\S^{3}$ different from the Hopf-link and let $\ul{q}=(q_{1},\dots,q_{n})$ be a family of integers;
  \begin{itemize}
    \item if $q_{i}=0$ for some $1\leq i \leq n$ then $L_{\ul{q}} = \S^{3}$;
    \item if $|q_{i}|\geq 2$ for all $1\leq i \leq n$ then there exists an irreducible representation
    \[
      \rho_{\ul{q}}:\pi_{1}L_{\ul{q}}\to \su.
    \]
  \end{itemize}
\end{proposition}

\begin{proof}
  First, if one of the $q_{i}$ is zero, the link $L_{(q_{1},\dots,q_{i})}$ is trivial so performing $1/{q_{k}}$-Dehn fillings for $i<k\leq n$ produces the standard $3$-sphere.

  On the other hand, if all the $|q_{i}|$ are greater than 1, by \Cref{prop:non_trivial_brunnian_links_fillings}, each $L_{(q_{1},\dots,q_{k})}$ for $k\leq n$ is nontrivial so $L_{(q_{1},\dots,q_{n-1})}$ is a nontrivial knot in $\S^{3}$ and Kronheimer-Mrowka's Theorem concludes the proof.
\end{proof}

By inclusion of $\su$ in $\sl$, we can consider $\rho_{\ul{q}}$ as an irreducible representation of $R(L_{\ul{q}})$ (with no nontrivial parabolic image).
Moreover, composing with the group homomorphism $\pi_{1}L\to\pi_{1}L_{\ul{q}}$, $\rho_{\ul{q}}$ may also be considered as an irreducible representation of $R(L)$.
The irreducible characters $\chi_{\ul{q}} = t(\rho_{\ul{q}})$ obtained this way are called \emph{Kronheimer-Mrowka characters} and we denote by $X_{\KM}(L)$ the Zariski closure in $X(L)$ of all Kronheimer-Mrowka characters:
\[
  X_{\KM}(L)=\ol{\left\{ \chi_{\ul{q}},~\ul{q}\in(\Z\setminus\{-1,0,1\})^{|L|} \right\}}.
\]
\begin{remark}
\label{rem:kro_mro_chars_several_comps}
  The subset $X_{\KM}(L)$ of $X(L)$ may contain several algebraic components.
\end{remark}
\begin{remark}
\label{rem:X_KM_natural}
  For any nontrivial, non-Hopf, Brunnian link $L=K\sqcup L_{0}$, the group homorphism $i_{q}:\pi_{1}L\to \pi_{1}L_{q}$ induces an algebraic map
  \[
    {i_{q}}^{\ast}: X(L_{q})\to X(L)
  \]
  and if $|q|\geq 2$, ${i_{q}}^{\ast}X_{\KM}(L_{q}) \subset X_{\KM}(L)$.
\end{remark}
Any representation $\rho_{\ul{q}}$ satisfies the $1/q_{K}$-Dehn filling relations for each component $K$ of $L$.
On the other hand, no $\rho_{\ul{q}}(\mu_{K})$ is trivial, since, otherwise, it would satisfy the $1/0$ relation on $K$;
it would then factor as a representation of $\S^{3}$ and therefore be trivial.
Since $\rho_{\ul{q}}$ factors in $\su$ this is equivalent to $\tr~\rho_{\ul{q}}(\mu_{K}\lambda_{K}^{q_{K}})=2$ and $\tr~\rho_{\ul{q}}(\mu_{K})\neq 2$.

It follows that any Kronheimer-Mrowka character $\chi_{\ul{q}}$ satisfies for any $K\subset L$:
\begin{align}
  I_{\mu_{K}\lambda_{K}^{q_{K}}}(\chi_{\ul{q}}) &= 2;\label{eq:1}\\
  I_{\mu_{K}}(\chi_{\ul{q}}) &\neq 2.\label{eq:2}
\end{align}

Finally, following \Cref{sec:eigenvalue_varieties_of_links_in_s_3}, we denote by $\evarGen{}{\KM}{L}$ the part corresponding to $X_{\KM}(L)$ in $\evar{L}$.
For any $\xi_{\ul{q}} \in \evarGen{}{\KM}{L}$ corresponding to a Kronheimer-Mrowka character $\chi_{\ul{q}}$ in $X_{\KM}(L)$, and any component $K$ of $L$, \eqref{eq:1} and \eqref{eq:2} imply that
\begin{align}
  \mf_{K}\lf_{K}^{q_{K}}(\xi_{\ul{q}}) &= 1;\label{eq:3}\\
  \mf_{K}(\xi_{\ul{q}})&\neq 1.\label{eq:4}
\end{align}

\begin{remark}
\label{rem:E_KM_non_reds}
  Together with the equations for $\aidGen{\red}{}{L}$ this implies that no such point $\xi_{\ul{q}}$ is in $\evarGen{\red}{}{L}$ so no component of $\evarGen{}{\KM}{L}$ is contained in $\evarGen{\red}{}{L}$.
\end{remark}



  \subsection{Maximal components} 
  \label{sub:maximal_components}

In this last section, we prove the following result which implies \Cref{th:brunnian_links_maximal_evar}:

\begin{theorem}
\label{th:kmc_maximal}
  For any nontrivial Brunnian link $L$ different from the Hopf link, $\evarGen{}{\KM}{L}$ contains a maximal component.
\end{theorem}

\begin{proof}
  This is proved by induction on the number of components of $L$.

  \noindent {\bf Base case: $L=K$}\\
  \indent
  For the base case, $L$ is a knot $K$ and the proof is the same as the one for the nontriviality of the $A$-polynomial of nontrivial knots from Dunfield-Garoufalidis in \cite{Dunfield:2004ei} or Boyer-Zhang in  \cite{Boyer:2005hu}.

  For any $|q|\geq 2$, performing $1/q$-surgery produces an irreducible character $\chi_{q}$ in $X(K)$ and a point $\xi_{q}=(m_{q},\l_{q})$ in $\evar{K}$.
  They show that there are infinitely many distinct $\l_{q}$ obtained this way so $\evarGen{}{\KM}{K}$ contains a curve different from the line $\lf-1$.
  We do not reproduce this proof here but very similar ideas are used for the induction step.

  \noindent {\bf Induction step: $L=K\sqcup L_{0}$}\\
  \indent
  Let $L=K\sqcup L_{0}$ be a nontrivial, non-Hopf, Brunnian link in $\S^{3}$.
  For any $|q|\geq 2$, $L_{q}$ is nontrivial, non-Hopf, and Brunnian, so we can assume, by induction, that $\evarGen{}{\KM}{L_{q}}$ contains a maximal component.

  We have the following commutative diagram:
  \[
    \xymatrix{
      X_{\KM}(L_{q}) \ar[r]\ar[d]_-{r_{q}}& X_{\KM}(L)\ar[d]^-{r}\\
      \prod_{J\neq K}X(T_{J}) & \prod_{J\subset L}X(T_{J})\ar[l]\\
      \prod_{J\neq K}E(T_{J})\ar[u]^-{d} & \prod_{J\subset L}E(T_{J})\ar[l]\ar[u]_-{d}
    }
  \]
  so there exists $X_{q}$ in $X_{\KM}(L)$ corresponding to $\evarC_{q}$ in $\evarGen{}{\KM}{L}$ such that $\dim~\evarC_{q}\geq |L|-1$.
  If $\dim~\evarC_{q} = |L|$ for some $q$ then there is nothing more to prove.

  Let us assume now that all the components $\evarC_{q}$ have dimension $|L|-1$.
  We will show that $\evarGen{}{\KM}{L}$ contains infinitely many different such subspaces $\evarC_{q}$;
  by algebraicity, this means that $\evarGen{}{\KM}{L}$ must contain a component of dimension $|L|$,= which will conclude the proof of \Cref{th:kmc_maximal}.

  The subspaces $\evarC_{q}$ will be separated using the following lemma:
  \begin{lemma}
  \label{lem:separate_E_q}
    For any integers $q,q'$,
    \[
      \evarC_{q}\subset \evarC_{q'}\Rightarrow {{\lf_{K}}^{q-q'}}_{|_{\evarC_{q}}}\equiv 1.
    \]
    Moreover, for any $|q|\geq 2$, the set $\{ p\in\Z ~|~{{\lf_{K}}^{p}}_{|_{\evarC_{q}}}\equiv 1\}$ is an ideal $d_{q}\Z$ with $q \not\in d_{q}\Z$.
  \end{lemma}
  \begin{proof}
    For any $\xi$ in $\evarC_{q}$, $\mf_{K}{\lf_{K}}^{q}(\xi)=1$ by equation \eqref{eq:3} so if $\xi$ also belongs to $\evarC_{q'}$, $\mf_{K}{\lf_{K}}^{q'}(\xi)=1$ and ${\lf_{K}}^{q-q'}(\xi) = 1$.
    Therefore, if $\evarC_{q}\subset \evarC_{q'}$, then ${\lf_{K}}^{q-q'} \equiv 1$ on $\evarC_{q}$.

    If $q$ is in the ideal $d_{q}\Z$, the surgery relation implies that ${\mf_{K}}_{|_{\evarC_{q}}} \equiv 1$, in contradiction with \eqref{eq:4}.
  \end{proof}

  If $S=\{q\in\Z\setminus\{-1,0,1\}~|~d_{q}=0\}$ is infinite then, by \Cref{lem:separate_E_q}, $\evarC_{q}\neq \evarC_{q'}$ for $q\neq q'$ in $S$, so $(\evarC_{q})_{q\in S}$ is a family of infinitely many distinct subspaces.

  Otherwise, there exists $N$ in $\N$ such that, for any $q\geq N$, $d_{q}\geq 2$.
  Let $(q_{i})_{i\in\N}$ be a family of integers such that:
  \begin{itemize}
    \item $q_{0} \geq N$;
    \item for any $j$ in $\N$, $q_{j+1}\geq q_{j}$ and $q_{j+1}\in\bigcap_{i=1}^{j} d_{q_{j}}\Z$.
  \end{itemize}
  Then, the following fact proves that $(\evarC_{q_{i}})_{i\in\N}$ contains infinitely many different subspaces:
  \[
    \forall~  i<j,~\evarC_{q_{i}} \neq \evarC_{q_{j}}.
  \]
  Indeed, for any $j$ in $\N$, let us assume that $\evarC_{q_{i}} = \evarC_{q_{j}}$ for some $i<j$.
  By \Cref{lem:separate_E_q}, this would imply that $q_{j}-q_{i}\in d_{q_{i}}\Z$;
  by construction, $q_{j}\in d_{q_{i}}\Z$ so this would imply $q_{i}\in d_{q_{i}}\Z$, a contradiction.

  We have proved that $\evarGen{}{\KM}{L}$ contains infinitely many different subsets of dimension $|L|-1$;
  by algebraicity, it must contain a component of dimension $|L|$, which concludes the proof of \Cref{th:kmc_maximal}.
\end{proof}






\renewcommand{\abstractname}{Acknowledgements}
\begin{abstract}
 The content of this paper forms part of the author's PhD thesis \cite{malabre_eigenvalue_2015}.
 He would like to express his gratitude to his advisors Michel Boileau and Joan Porti for their constant support during the realization of his PhD, and also thank Stephan Tillmann and Julien Marché for their valuable inputs on character and eigenvalue varieties.
\end{abstract}


\bibliographystyle{plain}
\bibliography{bib/brunnian}
\end{document}